\documentclass[journal,onecolumn]{IEEEtran}
\usepackage[lined]{algorithm2e}
\usepackage{xfrac} 
\usepackage{verbatim}
\usepackage{color}
\usepackage{multirow}
\usepackage{multicol}
\usepackage{qting}
\usepackage[normalem]{ulem}
\usepackage{enumitem}
\usepackage{graphicx} 
\usepackage{booktabs}
\usepackage[noadjust]{cite}
\usepackage{blindtext, graphicx,color, amsthm}
\usepackage[utf8]{inputenc}
\usepackage[english]{babel}
\usepackage{lettrine}
\usepackage{float}
\usepackage{epstopdf}
\usepackage{todonotes}
\usepackage[cmex10]{amsmath}
\usepackage{amsfonts}
\usepackage{mathtools}

\theoremstyle{definition}
\newtheorem{rmk}{Remark}

\newcommand{\correction}[1]{\textcolor{black}{#1}}

\def\undertilde#1{\mathord{\vtop{\ialign{##\crcr
				$\hfil\displaystyle{#1}\hfil$\crcr\noalign{\kern1.5pt\nointerlineskip}
				$\hfil\tilde{}\hfil$\crcr\noalign{\kern1.5pt}}}}}

\newcommand{\bd}{\boldsymbol}
\newcommand{\ie}{i.e.}
\newcommand{\eg}{e.g.}

\newenvironment{correction_E}{\par\color{black}}{\par}

\DeclareMathOperator*{\argmin}{arg\,min}
\makeatletter
\newcommand{\removelatexerror}{\let\@latex@error\@gobble}
\makeatother

\begin{document}
	
	\title{Tractable and Robust Modeling of Building Flexibility Using
		Coarse Data}

	\author{Jesus E. Contreras-Oca\~na,~\IEEEmembership{Student Member,~IEEE,} Miguel A. Ortega-Vazquez,~\IEEEmembership{Senior Member,~IEEE,} Daniel Kirschen,~\IEEEmembership{Fellow,~IEEE,}
		and~Baosen Zhang,~\IEEEmembership{Member,~IEEE}
		\thanks{The authors are with the Department of Electrical Engineering at the University of Washington (emails: \{jcontrer, kirschen, zhangbao\}@uw.edu) and with the Electric Power Research Institute (email: mortegavazquez@epri.com). This work is partially supported by the University of Washington Clean Energy Institute.}
	}

	\maketitle
	
	\begin{abstract}
		Controllable building loads have the potential to increase the flexibility of power systems. A key step in developing effective and attainable load control policies is modeling the set of feasible building load profiles.  In this paper, we consider buildings whose source of flexibility is their HVAC loads. We propose a data-driven method to empirically estimate a robust feasible region of the load using coarse data, that is, using only total building load and average indoor temperatures. The proposed method uses easy-to-gather coarse data and can be adapted to buildings of any type. The resulting feasible region model is robust to temperature prediction errors and is described by linear constraints. The mathematical simplicity of these constraints makes the proposed model adaptable to many power system applications, for example, economic dispatch, and optimal power flow. We validate our model using data from EnergyPlus and demonstrate its usefulness through a case study in which flexible building loads are used to balance errors of wind power forecasts.
	\end{abstract}

\begin{IEEEkeywords}
Buildings, flexibility, data-driven modeling. 
\end{IEEEkeywords}
	
	\IEEEpeerreviewmaketitle

		\begin{correction_E}
			
			\section*{Nomenclature}
			\addcontentsline{toc}{section}{Nomenclature}
\subsection*{Sets} \vspace{0pt}
\begin{IEEEdescription}[\IEEEusemathlabelsep\IEEEsetlabelwidth{$xxxxx$}]
	\item[$C$] Number of data clusters, indexed by $c$. 	
	\item[$ \bd{\mathcal{D}}$] Training dataset.
	\item[$ \bd{\mathcal{K}}$] Set of $K$ training days, indexed by $k$.
		\item[$ N$] Number of buildings, indexed by $i$.
	\item[$ \bd{\mathcal{P}}$] Feasible region of a building's load.
	\item[$T$] Number of time periods, indexed by $t$.
		\item[$\bd \Psi$] Set of explanatory variable data, its $k^\mathrm{th}$ element is $\psi_k$.
			 \item[$\bd \Omega^\mathrm{b}$] Set of load uncertainty scenarios, indexed by $\omega^\mathrm{b}$.
		\item[$\bd \Omega^\mathrm{w}$] Set of wind generation scenarios, indexed by $\omega^\mathrm{w}$.
\end{IEEEdescription}

\subsection*{Functions}
\begin{IEEEdescription}[\IEEEusemathlabelsep\IEEEsetlabelwidth{$xxxxx$}]
	\item[$ \bd g$] State transition function.
	\item[$\mathbb I$] Indicator function. 
	\item[$\bd \theta$] Zone temperatures.
	\item[$\hat \theta^\mathrm{L}$] Lower temperature estimate. 
	\item[$\hat \theta^\mathrm{U}$] Upper temperature estimate. 
	\item[$\xi$] Mapping of explanatory variables $\psi$ to its associated set of parameters $\Phi$.
\end{IEEEdescription}

\subsection*{Parameters and variables}
\begin{IEEEdescription}[\IEEEusemathlabelsep\IEEEsetlabelwidth{$xxxxx$}]
	 \item[$A$,$B$,$D$] Parameters of the RC circuit model. 
		 \item[$\underline{a}, \underline{b}$] Parameters of $\hat \theta^\mathrm{L}$. 
	   \item[$\overline{a}, \overline{b}$] Parameters of $\hat \theta^\mathrm{U}$.
	   \item[$J^\mathrm{A}$]  Measure of tightness of the prediction band. 
	   \item[$J^\mathrm{U}/J^\mathrm{L}$] Upper/lower estimate error.
	   \item[$M$] Arbitrary integer greater than $1$.
	    \item[$\bd p$]  Load vector, its $t^\mathrm{th}$ element is $p_t$.
	    \item[$p^\mathrm{base}$] Base load.
	    \item[$p^\mathrm{hvac}$] HVAC load. 
	     \item[$\bd u$] Thermal input.
	     \item[$v$] Compensation for wind power balancing.
	       \item[$\bd W$] Matrix of normalized training data for clustering, its $k^\mathrm{th}$ column is $\bd w_k$
	     \item[$\bd x$] Thermal state of the building.
	     \item[$\alpha$] Robustness parameter. 
	      \item[$\bd \beta$] $M$-sized vector whose $i^\mathrm{th}$ entry, $\beta_i$, is $\frac{i-1}{M-1}$.
	     \item[$\bd \Delta$] Wind forecast error.
	     \item[$\bd \epsilon$] Stochastic component of the load. 
	     \item[$\bd \nu$] Wind generation. 
	     \item[$\pi^\mathrm{out}$] Portion of measurements outside the prediction band.
	     \item[$\bd \tau$] Vector of energy prices, its $t^\mathrm{th}$ element is $\tau_t$.
	     \item[$\bd{\Phi}$] 6-tuple that contains $\left( \hat{   p}^\mathrm{min} , \hat{p}^\mathrm{max},   \hat{ \theta}^\mathrm{min},  \hat{ \theta}^\mathrm{max}, \hat \theta^\mathrm{L}, \hat \theta^\mathrm{U} \right)$. 
	     \item [$\bd \phi^\mathrm{in}$] Average indoor temperature, its $t^\mathrm{th}$ element is $\phi^\mathrm{in}_t$.
	     \item [$\bd \phi^\mathrm{out}$] Outdoor temperature, its $t^\mathrm{th}$ element is $\phi^\mathrm{out}_t$ .
	     \item [$\bd  \Sigma$] Covariance matrix of $\bd \epsilon$.
\end{IEEEdescription}

\subsection*{Accents and subscripts} 
\begin{IEEEdescription}[\IEEEusemathlabelsep\IEEEsetlabelwidth{$xxxxx$}]
\item[$\widehat{x}$] Approximation of $x$.
\item[$x_{a:b}$ ] Vector of composed of the $a^\mathrm{th}$ through the $b^\mathrm{th}$ elements of vector $x$.
\item[$x^\mathrm{max}$ ] Upper limit of $x$.
\item[$x^\mathrm{min}$ ] Lower limit of $x$.
\item[$x^\top$ ] Transpose of $x$.
\end{IEEEdescription}
			
		\end{correction_E}

	\section{Introduction}

\begin{correction_E}\IEEEPARstart{P}{ower} system flexibility is defined as the ability to respond to changes in demand or supply within a given time frame~\cite{cochran2014flexibility}.  Traditionally, the main (and usually sole) source of flexibility in power systems has been flexible generation resources, \eg, simple and combined cycle gas turbines. Meanwhile, the load has been treated as a fixed quantity to be followed by the flexible supply-side. Nowadays, however, thermostatically controlled loads (TCL) in buildings (\eg, HVAC units, refrigerators, and water heaters) have emerged as important sources of flexibility~\cite{LOPES20161053,Ma2012}.  In this new environment, TLCs can provide flexibility from the demand-side by altering their consumption to accommodate power variations. \end{correction_E}

    \begin{correction_E} Some of the benefits of increased power system flexibility include deferral of infrastructure investments~\cite{Albadi,contreras2017non}, increase of renewable energy hosting capacity~\cite{Lund_Lindgren_2015}, increased economic efficiency~\cite{GOY20153391}, and others~(see, e.g., \cite{7416149, bebic2015} and the references therein).  However, to fully harvest the flexibility of TCLs, challenges still remain: developing appropriate building models, aggregating those models for large-scale implementation, data privacy, state estimation of flexible loads, among others~\cite{bebic2015, Piette2009,DOE_benefits}. This paper focuses on overcoming these challenges by proposing a method that uses easy-to-collect data from buildings to find \textbf{tractable} and \textbf{robust building models} that capture their \textbf{load flexibility}.

The concepts in bold font in the previous paragraph seem simple but are loaded with meaning.  We start by introducing the concept of load flexibility. \correction{Similar to~\cite{Muller2015}, we define \textbf{load flexibility}, or equivalently, the feasible region of the load $\bd{ \mathcal{P}}$, as the collection of load profiles that satisfy the user requirements, \eg, thermal comfort, technical limits.}

We then turn to the concept of robustness. A model of $\bd{ \mathcal{P}}$,  denoted by $\widehat{\bd{ \mathcal{P}}}$, is said to be \textbf{robust} if an arbitrary element (\ie, a load profile) of $\widehat{\bd{ \mathcal{P}}}$ is also contained in $\bd{ \mathcal{P}}$ to a degree of certainty. This feature is of particular importance since it ensures, to said degree of certainty, that a load profile in the model is actually attainable by the physical building.

In this work, a model $\widehat{\bd{ \mathcal{P}}}$ is said to be  \textbf{tractable} if it can be easily incorporated into a desired power system analysis frameworks.  For instance, since the unit commitment problem is typically modeled as a mixed-integer linear program (MILP), to seamlessly incorporate flexible loads, the feasible set of loads should be described by linear constraints. Other typical power system analysis frameworks include optimal power flow, economic dispatch, etc.~\cite{cochran2014flexibility, Wang_2003_Demand, Wang_2017_Fully}. 

\end{correction_E}
\subsection{Flexibility of HVAC loads}
We consider buildings whose source of flexibility is their HVAC loads\footnote{We consider heating and cooling loads because they are the largest component of commercial loads~\cite{CBECS_2012}.} and indoor temperature as the only controllable comfort index. Since indoor temperature is typically allowed to exist within an admissible range for human comfort,~\eg, from $20^\circ\mathrm{C}$ to $25^\circ\mathrm{C}$, there exists a range of HVAC load levels that achieve admissible indoor temperatures.  Then, the building operator could choose a load level among the set of possibilities that accomplishes a power system-level objective, \eg, demand response, while satisfying indoor temperature requirements. However, the relationship between indoor temperature and the electrical power consumption can be complex~\cite{Bacher2011}, and finding the set of load profiles that maps to admissible temperatures is nontrivial.

This paper presents a \emph{data-driven} approach to \emph{appropriately} model the feasible region of a building load. We are interested in developing models of demand-side flexibility that \emph{i)} are computationally tractable, \emph{ii)} do not compromise occupant comfort, and \emph{iii)} can be estimated with relatively simple and easy-to-obtain sets of data. The proposed model is described by a set of linear constraints and continuous variables, making it easy to integrate into existing power system analysis frameworks. It is also robust in the sense that it ensures  (to a certain degree of confidence) that a load profile is attainable without violating indoor temperatures limits. Finally, this approach uses small amounts of relatively coarse\footnote{Gathering fine-grained and/or architectural building parameters may be expensive or not viable.}, high-level data: average indoor temperatures rather than zone-specific temperatures and total building load rather than individual appliance loads. The characteristics of the required data make our approach attractive to entities that manage buildings that are not instrumented at the appliance level (which is still the case for most buildings).  As a consequence of its data-driven nature, our approach does not rely on human inputs of architectural building parameters.

\subsection{Related works}
Perhaps the most popular way of modeling building thermal dynamics, \ie, the relation between indoor temperature and HVAC load, is the \emph{resistance-capacitance (RC) circuit model}~\cite{Ma2012,Radecki,Gouda2002,Hao2015}. Typically, the values of the resistance and capacitance parameters are calculated from the building's specifications, including zonal volume, insulation material, and wall area, among others. These models provide an easy-to-use linear representation of a building's thermal dynamics. However, they are restrictive and calculating the parameters may be costly and labor intensive. In contrast, our approach uses building-level metered data to identify a thermal model that is just as easy to use as the RC models and do not require human input.

Using metered data to identify building models have been studied~\cite{  Bacher2011, Hughes2016,Goyal2011,Gunay2017, Royer2014}. The authors of~\cite{Bacher2011} propose a maximum likelihood estimation of the RC parameters of a building. In~\cite{Hughes2016}, the authors propose a method to identify the parameters of a building as a ``virtual battery'' for system-wide frequency regulation.  Both models in~\cite{Bacher2011} and~\cite{Hughes2016} are suited for short prediction horizons, \eg, 5 minutes. However, we are interested in longer prediction horizons for power system operation purposes, \eg, 24 hours. The authors of~\cite{Goyal2011} propose an estimation of intrabuilding RC parameters to predict the internal thermal dynamics of a building using room-level temperature data. Here, we are interested in models that use coarser data (\eg, average zonal temperature rather than individual room temperature) that are easier to obtain, especially for older buildings.

In addition to data efficiency, we are interested in developing mathematically \emph{simple} models that can be easily adopted in a wide range of power system applications. Thus, unlike the artificial neural network models in~\cite{Gunay2017,Chen_2017_modeling} that employ nonlinear activation functions, our proposed model employs exclusively linear equalities and inequalities that can be easily embedded in most power system optimization and control settings. The work in \cite{Royer2014} is closely related to ours since it also identifies a coarse (\eg, facility-level) models of thermal dynamics. The main difference is that the approach in \cite{Royer2014} provides a single central temperature estimate, whereas ours provides upper and lower estimates.

\subsection{Overview of the proposed method}
\begin{figure}[tb]
    \centering
    \includegraphics[width=0.35\textwidth]{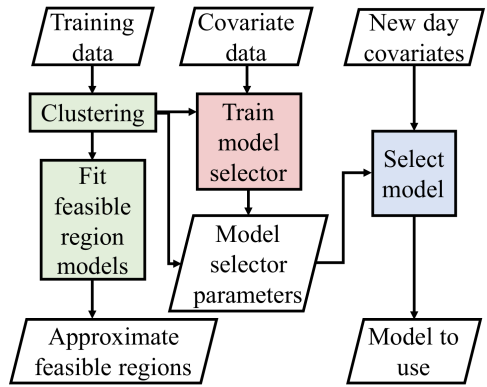} 
    \caption{High-level illustration of the proposed method. First (green blocks), we group similar training data points into a number of clusters. Each cluster is used to fit a model of the feasible region. Then (red block), the explanatory variables and the training data clusters are used to train the parameter of ``model selecting function.''  Finally (blue block)  we feed expected explanatory variables of a new day to select the feasible region model to use. Note that the first two steps (green and red blocks) are performed when fitting the models while the last step is performed when using the model to predict the flexibility of the building.}  \label{fig:high_level_alg}
\end{figure}

The method to find appropriate models of the feasible region of a building load is comprised three major tasks: clustering the training data, fitting the feasible region models, and selecting the model to use. First, we group similar training data points into a number of clusters. Each cluster is then used to fit a model of the feasible region of the building's load. Each training data point is associated with a set of explanatory variables, \eg, outdoor temperature, solar irradiation, or day of the week. The explanatory variables and the training data clusters are then used to train a ``model selecting function''. The model selecting function takes the explanatory variables associated with a training data point and predicts the best model to use.  Finally, we feed expected explanatory variables of a new day, \ie, a day outside the training set, to select the feasible region model to use (for such new day). Fig.~\ref{fig:high_level_alg} illustrates these three major tasks.

The major contributions of our work are:
\begin{itemize}
    \item A method to describe the feasible region of a building's load that uses a small amount of easy-to-obtain data. We first group the training data into $C$ clusters of similar days. Clustering allows us to segregate the training data by classes of thermal behaviors and train models for each cluster rather than a single general model. Then, we train a linear but robust model of building thermal dynamics using a technique we call bounded least squares estimation~(BLSE). Rather than providing a central prediction of the building indoor temperature, the BLSE provides a prediction band, \ie, upper and lower estimates of indoor temperature.
    \item Validation of our models using data from the building modeling software EnergyPlus~\cite{crawley2001}.
    \item A demonstration showing how the proposed model can be used to schedule building loads to mitigate the discrepancies between expected and actual wind power generation.
\end{itemize}

\subsection{Organization of this paper}
The rest of this paper is organized as follows. Section~\ref{sec:underlying_bldg} 
describes the model of a generic building and defines the feasible region of the load and Section~\ref{sec:TRAFL}
 introduces the tractable and robust model of the feasible region.  Section~\ref{sec:data} describes the data used and outlines the procedure for estimating the approximation of the feasible region. Section~\ref{sec:case} validates the model and compares it with the traditional RC circuit model. Section~\ref{sec:wind_power_bal} shows how the proposed approximation can be used to model a building that uses its flexibility to compensate wind power forecasting errors. Section~\ref{sec:conclsion} concludes this paper.

	\section{Preliminaries}
	\label{sec:underlying_bldg}
We define the feasible region of the load as the set of all load profiles that meet power and indoor temperature limits. The maximum power limit is the non-HVAC building load (or base load) plus the installed capacity of the HVAC system. The minimum power limit, on the other hand, is the base load plus the minimum power of the HVAC system. The temperature limits are given by predefined comfort limits.

\subsection{Building thermal dynamics} \label{subsec:thermal_dyn}
 The building's thermal dynamics model describes the behavior of the indoor temperatures of a building as a function of the heating,  cooling, and internal and external disturbances. For each thermal zone of the building, there are two quantities of interest: the stored energy and the temperature. The stored energy, denoted by $\bd x$, represents the thermal state of the building. The temperature of each zone is denoted by $\bd \theta$ and is constrained by the comfort range of the users in the building. The input to the system is denoted by $\bd u$, representing the thermal inputs of the HVAC (energy is injected when heating and withdrawn when cooling). Then the states $\bd x$ evolves as:
\begin{equation*} \label{eq:state_transition}
    \bd x_{t} = \bd g_{t}(\bd x_{t-1},\bd u_{t-1}),
\end{equation*}
where $\bd g$ is the state transition function (not necessarily linear nor observable). The indoor temperature at time $t$ is then modeled as a function of the state and the input at the current time, which we write as $\bd \theta_{t}(\bd x_t, \bd u_t)$.

The thermal input $\bd u_t$ is not a directly controllable in practice. Instead, we want to link the temperature to the electrical load of a building (in units of kW). Since we assume that the HVAC is the only controllable load, we can express the total building load, $p_{t}$,  as the sum of the HVAC load and the base load: $ p_{t}= p^\mathrm{hvac}_{t}  + p^\mathrm{base}_{t}. $

The thermal input of the HVAC system, $\bd u_t$, is a function of the electrical load $p_t$. Then, we can express the indoor temperature at time $t$ as a function of the load at $t$, the load at $t-1$ and the state at $t-1$, written as
$\bd \theta_{t}(\bd x_{t-1}, p_{t-1}, p_{t})$. Rolling out this relationship backward in time until $t=0$, we can eliminate the dependencies on $\bd x_1,\dots,\bd x_T$, and write $\bd \theta_{t}$ as function of $\bd x_0$ and $p_0, p_1,\dots,p_t$,
\begin{equation*}
\bd \theta_{t}(\bd x_0, \boldsymbol{p}_{0:t}),
\end{equation*}
where $\boldsymbol{p}_{0:t} =\begin{bmatrix} p_0&\dots&p_t\end{bmatrix}^\top$.

\subsection{Feasible region of the load}
\begin{figure}
    \centering
    \includegraphics[width=0.45\textwidth]{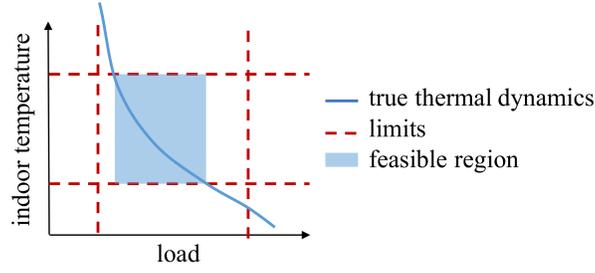} 
    \caption{Illustration of the feasible region of the load and a (hypothetical) thermal dynamics function. With a complicated non-linear dynamics function such as the one illustrated here, the feasible region described by Eqs.~\eqref{eq:feasible_region} 
 is too complex to be used in typical power system frameworks. } \label{fig:true_feasible_region_illustration}
\end{figure}
Since we assume that the building operation is constrained by both load and temperature limits, the feasible set of load profiles is given by
     \begin{subequations}
        \begin{align}
        \bd{\mathcal{P}} =\{   \bd{p}\; | \; \label{eq:feasible_region_1} &
         p^\mathrm{min}_t\le  p_t \le  p^\mathrm{max}_t \; \forall\; t=1, \dots, T\\
         & \bd \theta_{t}(\bd x_0 , \boldsymbol{p}_{0:t})  \ge  \bd{\theta}_t^\mathrm{min}  \; \forall\; t=1, \dots, T \\
         &  \bd \theta_{t}(\bd x_0, \boldsymbol{p}_{0:t}) \le    \bd{\theta}_t^\mathrm{max} \; \forall\; t=1, \dots, T 
          \label{eq:feasible_region_3}
            \\ 
         & \bd{p}_{0:t} = \begin{bmatrix} p_0&\dots&p_t\end{bmatrix}^\top, \; \bd p  = \bd p_{1:T} \} . \label{eq:feasible_region_4}
        \end{align}\label{eq:feasible_region}
         \end{subequations}
The symbol $\bd{\mathcal{P}}$ describes set of load profiles that are within the minimum and maximum load limits ($\bd p^\mathrm{min}$  and  $\bd p^\mathrm{max}$) and whose associated indoor temperatures are within admissible comfort limits ($\bd{\theta}_t^\mathrm{min}$ and $ \bd{\theta}_t^\mathrm{max} $). The region $\bd{\mathcal{P}}$ is hard to characterize mainly due to the complicated thermal dynamics function, $\bd{\theta}_t$. This motivates the main objective of this work: finding a good approximation $\bd{\mathcal{P}}$. Fig.~\ref{fig:true_feasible_region_illustration} illustrates the feasible region in one load dimension.

			\section{A tractable and robust approximation of the feasible region}
			\label{sec:TRAFL}
			We look for two important features of an approximation: robustness and simplicity. The former is important because, in the vast majority of applications, the main purpose of the HVAC system is to maintain acceptable levels of occupancy comfort, while servicing the grid is of secondary priority (if at all). The latter is important because, to interface with the larger electrical grid, the description of $\bd{\mathcal{P}}$ should be tractable in grid-related optimization and control problems. 

        \begin{figure}        
        \centering
        \includegraphics[width=0.485\textwidth]{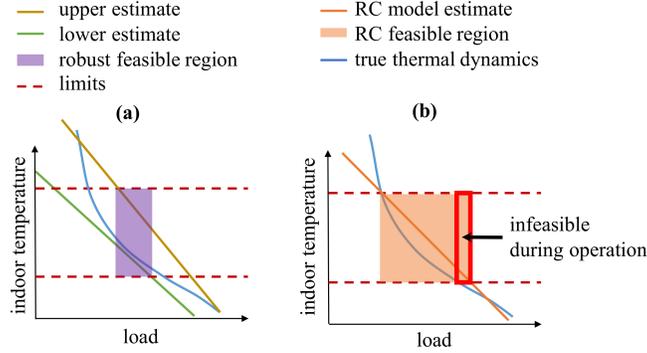} 
        \begin{correction_E}
            \caption{ Illustration of the proposed robust approximation of the feasible region of the load (plot a) and the feasible region of the load given by the RC circuit model (plot b). In the proposed approximation, the upper temperature estimate is upper bounded by the maximum temperature limit. Similarly, the lower temperature estimate is lower bounded by the minimum temperature limit. The RC circuit model, on the other hand, limits a central estimate to be within minimum and maximum temperatures. Note that the RC circuit model may overestimate the feasible region which may lead to load profiles that are infeasible during actual building operation.  }  \label{fig:model_vs_RC_illust}
        \end{correction_E}
        \end{figure}

We attain robustness by producing an indoor-temperature \emph{prediction band} of arbitrary confidence. Then, we limit an \emph{upper estimate} of the indoor temperature to be under the maximum temperature limit and a \emph{lower estimate} to be over the minimum limit.  This stands in contrast to the typical RC circuit model which produces a \emph{central estimate} of indoor temperature~\cite{Ma2012,Radecki,Gouda2002,Hao2015,Contreras-Ocana2016}.  Thus, when the temperature is underestimated in the RC circuit model, the maximum temperature limit could be violated. Similarly, the minimum temperature limit could be violated when the temperature is overestimated.

A natural implication of using a prediction band rather than a central estimate is that it allows us to handle a diverse set temperature dynamics functions. Rather than attempting to fit a linear function to potentially complicated underlying thermal dynamics, our method determines linear upper and lower estimates of the thermal dynamics function. Fig.~\ref{fig:model_vs_RC_illust} illustrates how a linear prediction band can handle a complicated thermal dynamics function.

A ``simple'' approximation should be tractable and require coarse-grained data. Tractability is important to seamlessly incorporate our model of flexibility into power system frameworks that exist in practice, e.g., the UC problem, optimal power flow, among others. For instance, the PJM Interconnection and California ISO implement their UC problems as MILPs~\cite{ott2010unit,Rothleder2010unit}. Similarly, most European market designs are implemented using MILPs~\cite{Chatzigiannis2016European}. Thus, we achieve tractability by approximating the feasible region of the load, $\bd{\mathcal{P}}$, as a polyhedron described by \emph{linear} relations. This way, a flexible load can be easily incorporated into the previously mentioned frameworks as a variable constrained by a  polyhedral approximation of $\bd{\mathcal{P}}$. Finally, our approximation is low-dimensional since we only use coarse data: average indoor temperatures and building-level load. 

    \begin{correction_E}
        \begin{rmk}
            In this work, zonal indoor temperatures are weighted by the zone's volume to determine the building's average indoor temperature.
        \end{rmk}
\end{correction_E}

The mathematical simplicity of our model stands in contrast to neural network-based models like the ones in~\cite{Gunay2017,Chen_2017_modeling}. While such models are useful for certain applications, \eg, local load control, their non-linear representations makes them ill-suited for UC and market models.

Let our robust approximation feasible region of a building's load be denoted by
\begin{subequations}\label{equation:aprox_feasible_reg}
    \begin{align} 
    \widehat{\bd{\mathcal{P}}}=\{  \bd{p} \; | \; & \hat{ p}^\mathrm{min}_t \le   p_t \le \hat{ p}^\mathrm{max}_t  \; \forall\; t=1, \dots, T\\
    & \hat{ \theta}^\mathrm{U}_{t}(\phi_0^\mathrm{in}, \phi_t^\mathrm{out}, \bd p_{1:t})  \le \hat{  \theta}_{t}^\mathrm{max} \; \forall \; t=1,\dots,T \\
    &\hat{\theta}^\mathrm{L}_{t}(\phi_0^\mathrm{in}, \phi_t^\mathrm{out}, \bd p_{1:t})  \ge \hat{\theta}_{t}^\mathrm{min}\; \forall \; t=1,\dots,T \\ 
    & \bd{p} = \begin{bmatrix} p_1&\dots&p_T\end{bmatrix}^\top, \; \bd p = \bd p_{1:T} \}.
    \end{align} 
\end{subequations} 
This approximation has a similar structure to the feasible region described by Eqs.~\eqref{eq:feasible_region}.   In this case, however, the load $\bd p$ is constrained by approximations of the load limits ($\hat{\bd p}^\mathrm{min}$ and $\hat{ \bd p}^\mathrm{max}$). A more notable difference is that indoor temperatures are described by upper and lower estimates, ($\hat{ \theta}_t^\mathrm{U}$ and $\hat{ \theta}_t^\mathrm{L}$) and limited by approximations of the maximum and minimum temperature limits ($\hat{\theta}^\mathrm{min}_t$ and $\hat{\theta}^\mathrm{max}_t$), respectively. Fig.~\ref{fig:model_vs_RC_illust}(a) illustrates the approximation of the feasible region in one load dimension.

To achieve the tractability property previously described, we model $\hat{\theta}_t^\mathrm{U}$ and $\hat{ \theta}_t^\mathrm{L}$ as \emph{affine} functions of the initial indoor temperature $\phi_0^\mathrm{in}$, outdoor temperature $\phi_t^\mathrm{out}$, and load from $1$ to $t$:
\begin{subequations}
    \begin{align} 
    &\hat \theta^\mathrm{U}_{t}(\phi_0^\mathrm{in},\phi_t^\mathrm{out}, \bd p_{1:t}) = \overline{\bd{a}}_{t}^{ \top}\bd{p}_{ 1:t}+\overline{\bd{b}}_{t}^{ \top}\begin{bmatrix}\phi_{0}^\mathrm{in} & \phi^\mathrm{out}_{t}  & 1\end{bmatrix}^\top \label{eq:u_l_estimate_functions_u}
    \\
    &\hat \theta^\mathrm{L}_{t}(\phi_0^\mathrm{in},\phi_t^\mathrm{out}, \bd p_{1:t}) = \underline{\bd{a}}_{t}^{ \top}\bd{p}_{ 1:t}+\underline{\bd{b}}_{t}^{ \top}\begin{bmatrix}\phi_{0}^\mathrm{in} & \phi^\mathrm{out}_{t}  & 1\end{bmatrix}^\top. \label{eq:u_l_estimate_functions_l} 
    \end{align} \label{eq:u_l_estimate_functions}
\end{subequations}
The vectors $\overline{\bd{a}}_{t}, \;\underline{\bd{a}}_{t} \in \mathbb{R}^{t}$ relate the building load from time $1$ to $t$, $\bd p_{1:t}$, to upper and lower estimates of indoor temperature at time $t$, respectively. The vectors $\overline{\bd{b}}_{t},\;\underline{\bd{b}}_{t} \in \mathbb{R}^{3}$, on the other hand, relate outside ambient temperature and the initial indoor temperature to upper and lower approximations indoor temperature at time $t$, respectively. The last element of $\overline{\bd{b}}_{t}$ and $\underline{\bd{b}}_{t}$ is the offset of their respective functions. 

For each time period, the feasible region is described by the 6-tuple $\bd \Phi_t = \left( \hat{   p}^\mathrm{min}_t , \hat{p}^\mathrm{max}_t,   \hat{ \theta}^\mathrm{min}_t,  \hat{ \theta}^\mathrm{max}_t, \hat \theta^\mathrm{L}_t, \hat \theta^\mathrm{U}_t \right)$ and the collection of all $\bd \Phi_t$'s from $t=1$ to $t=T$ describe the entire feasible region. The next section shows how to find $\bd \Phi_t$ from data.

	\section{Estimating the  feasible region}
	\label{sec:data}
	
\correction{We use time-series data of total building load, indoor temperature, and outdoor temperature to learn an approximation of $\bd{\mathcal{P}}$. The data is simulated on EnergyPlus\footnote{\correction{EnergyPlus is widely used in the literature in lieu of actual building measurements that are rarely available to academic researchers, \eg, in References~\cite{Bojic,Royer2014,yoon2014dynamic}}} and based on typical commercial buildings detailed in the report U.S. Department of Energy Commercial Reference Building Models of the National Building Stock~\cite{CommercialRef}. While the data itself is simulated, the underlying data to construct the models  is real and characteristic of buildings in the United States.}

\begin{figure}
 \centering
 \includegraphics[width=0.485\textwidth]{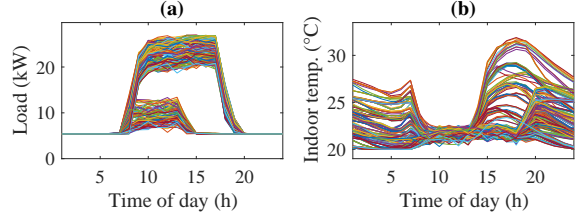} 
 \begin{correction_E}
 \caption{Load and indoor temperature time-series data for a small office building simulated in EnergyPlus. Plot a shows load data for 300 different summer days and plot b shows indoor temperature. There are two distinct types of days corresponding to weekday and weekends. Temperatures that might not comply with occupant comfort (either too high or too low) are observed during weekends and after business hours when occupants are not present in the building. }  \label{fig:training_data}
 \end{correction_E}
\end{figure}

Fig.~\ref{fig:training_data} shows load and indoor temperature data for a small office building for $300$ summer days, produced using EnergyPlus. \correction{In our work, we find that using $300$ training days is sufficient to identify most of the behaviors.  Since buildings are typically designed to operate for years without major renovations~\cite{Inaba2014}, lack of data should not be a problem after a few months of operation.}

While there are other valuable pieces of information that could help to better approximate $\bd{\mathcal{P}}$, we are interested in using a small amount of data. For instance, using HVAC rather than total load could result in a more accurate approximation of the thermal dynamics. However, such data may not be readily available or gathering it may be difficult and costly~\cite{Knight12}.

Denote the set of training data points for each set of parameters $\bd \Phi_t$ as $\bd{\mathcal{D}}_t=\left\{ \left(\boldsymbol{p}_{k,1:t} ,\phi^\mathrm{in}_{k,0}, \phi^\mathrm{in}_{k,t}, \phi^\mathrm{out}_{k,t}\right) \right\}_{k\in \bd{\mathcal{K}}}$.  Here $\bd{\mathcal{K}}$ is the set of $K$ training days and the subscript $k$ denotes data of the $k^\mathrm{th}$ day.

\subsection{Clustering}
Needless to say, each one of the training days is different to each other. For instance, the outdoor temperatures of two days are never \emph{exactly} the same. Therefore, each day is actually associated with a different region $\bd{\mathcal{P}}_k\;  \forall\; k \in \bd{\mathcal{K}}$. However, fitting one model per day is difficult and of little use since no future day is exactly like any previous day. Instead, we fit $C_t$ different values of the parameters $\bd \Phi_t$, with the goal of capturing distinct types of days. In general,  $C_t$ is much smaller than the total number of training days.

Denote the $C_t$ different parameters of the load-indoor temperature relation, load, and temperature limits as $\left\{  \bd \Phi_{1,t}, \bd \Phi_{2,t}, \dots, \bd \Phi_{C_t,t}\right \}$ where
\begin{equation*}
\bd \Phi_{c,t} = \left( \hat{   p}^\mathrm{min}_{c,t} , \hat{p}^\mathrm{max}_{c,t},   \hat{ \theta}^\mathrm{min}_{c,t},  \hat{ \theta}^\mathrm{max}_{c,t}, \hat \theta^\mathrm{L}_{c,t}, \hat \theta^\mathrm{U}_{c,t} \right).
\end{equation*}
 Intuitively, one would like set of parameters to model days that are similar to each other. We accomplish this by clustering days that exhibit similar load-indoor temperature-outdoor temperature relationships. Specifically, we use the K-means clustering algorithm to partition the training set $\bd{\mathcal{D}}_t$ into subsets used to train the different models. The training set contains data in different units (temperature and power units) and likely different magnitudes. To accommodate this, we normalize every dimension of the training set to have an $\ell_2$ norm of $1$ before applying the K-means algorithm. We denote the resulting clusters of $\bd{\mathcal{D}}_t$ as 
 $\{\bd{\mathcal{D}}_{1,t}, \bd{\mathcal{D}}_{2,t},\hdots,\bd{\mathcal{D}}_{C_t,t}\}$
 where each set $\bd{\mathcal{D}}_{c,t}$ is associated with a subset of $\bd{\mathcal{K}}$ denoted as $\bd{\mathcal{K}}_{c,t}$. Then the functions $\left(\hat{\theta}_{c,t}^\mathrm{L},\hat{\theta}_{c,t}^\mathrm{U} \right)$  and the rest of the parameters in $\bd \Phi_t$ are trained from the data in $\bd{\mathcal{D}}_{c,t}$. The clustering method is detailed in Appendix~\ref{app:data}.

 \begin{figure}
    \centering
    \includegraphics[width=0.30\textwidth]{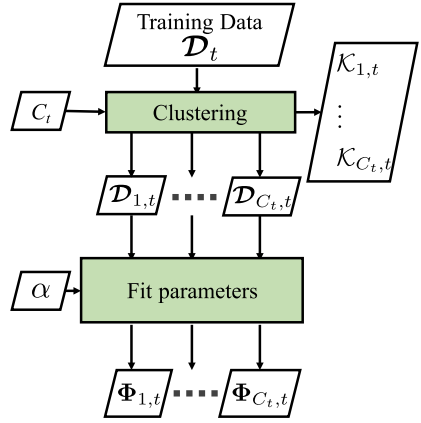}
    \caption{Illustration of the algorithm that groups the training dataset $\bd{\mathcal{D}}_t$ into $C_t$ clusters and trains the parameters $\{\bd{\Phi}_{1,t}, \bd{\Phi}_{2,t},\hdots,\bd{\Phi}_{C_t,t}\}.$  Each subset of training data $\bd{\mathcal{D}}_{c,t}$ is used to train its corresponding set of feasible region parameters $\bd{\Phi}_{c,t}$.}  \label{fig:algorithm}
\end{figure}

\subsection{Robust model of the load-indoor temperature relationship} \label{subsec:hypothesis}

The inputs of this portion of the algorithm are the training data $\bd{\mathcal{D}}_t$, the number of clusters $C_t$, and a robustness tuning parameter $\alpha \in (0,1)$.  \correction{The robustness parameter $\alpha$ represents the proportion of temperature observations outside the upper and lower predictions. Thus, a smaller alpha leads to wider, more robust prediction bounds. Conversely, a larger alpha leads to tighter prediction bounds. The Case Study demonstrates how a small $\alpha$ leads on the one hand, to more aggressive provision of flexibility but on the other hand, a higher risk of violation of the temperature limits. Appendix~\ref{app:blse} details the mechanism whereby $\alpha$ influences the robustness of the prediction bounds. }

\begin{figure}
    \centering
    \includegraphics[width=0.4\textwidth]{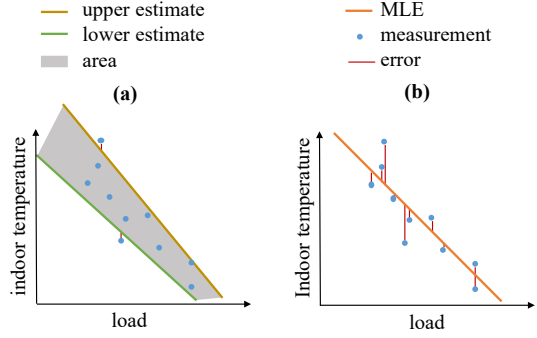}
    \caption{Illustration of (a) the proposed bounded least-squares estimation and (a) classical least-squares estimation.  The former provides an upper and a lower prediction that minimize two weighted objectives: 1) the area between the predictions and 2) the MSE of the points outside the prediction bounds. The latter provides a central maximum likelihood (MLSE) linear estimator.}  \label{fig:BLSE_est}
\end{figure}

\begin{figure}
    \centering
    \includegraphics[width=0.485\textwidth]{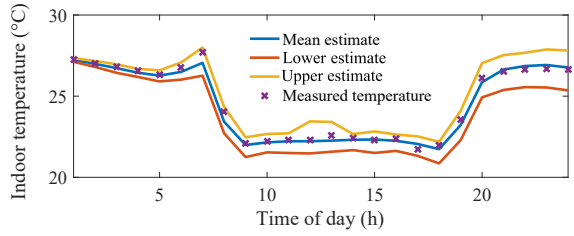} 
    \caption{Measured indoor temperature, indoor mean temperature estimate, and prediction bounds for a sample day.}  \label{fig:sample_day_prediction}
\end{figure}

We learn each set of parameters $\bd \Phi_{c,t}$ using its respective training data subset $\bd{\mathcal{D}}_{c,t}$ as illustrated in Fig.~\ref{fig:algorithm}. The parameters of the load -indoor temperature relation, $\overline{\bd{a}}_{c,t}$, $\underline{\bd{a}}_{c,t}$, $\overline{\bd{b}}_{c,t}$, and $\underline{\bd{b}}_{c,t}$, are trained using a least squares estimation (LSE)-inspired algorithm that we call ``bounded least squares estimation (BLSE)." The classic LSE calculates a line that minimizes the mean squared error (MSE) of the prediction. The BLSE finds \emph{two} lines (an upper and a lower prediction) such that a weighted sum of two objectives is minimized: \emph{1)} the squared error of the points outside of the bounds and \emph{2)} a measure of the area between the predictions. The weights assigned to each objective, determined by the robustness tuning parameter $\alpha$, influences the tightness of the bounds: the higher the weight assigned to the area objective is, the tighter the predictions are. The tightness of the predictions affects the load scheduling problem: overly tight predictions may translate into overestimation of the building flexibility while looser predictions might translate into an overly conservative approximate feasible region of the load. Figure~\ref{fig:BLSE_est} shows an illustration of classical LSE and contrasts it with the proposed BLSE. Figure~\ref{fig:sample_day_prediction} shows measured indoor temperature, mean indoor temperature estimate, and prediction bounds for a sample day. The BLSE algorithm is detailed in Appendix~\ref{app:blse}.
.

\subsection{Estimates of the temperature and load limits}
\label{subsec:limits}
In some ways, estimating the parameters of the functions $\hat  \theta^\mathrm{U}_{c,t}$ and $\hat \theta^\mathrm{L}_{c,t}$ is easier than estimating the temperature and load limits.  The former is a supervised learning problem, while the latter is unsupervised since we do not directly observe the limits. Therefore, we approximate the temperature and load limits of cluster $c$ as the highest/lowest observed values during days in $\bd{\mathcal K}_{c,t}$. Appendix~\ref{app:est} offers details about this method.

\begin{correction_E}
\subsection{Model selection}

 Suppose we would like to estimate the building's flexibility during day $K+1$, \ie, during a day outside the training set. The first question is: for each time period $t$, which set of parameters in $\left\{  \bd \Phi_{1,t}, \bd \Phi_{2,t}, \dots, \bd \Phi_{C_t,t}\right \}$  should make up the approximate region $\widehat{\bd{\mathcal{P}}}_{K+1}$? 

\begin{figure}
    \centering
    \includegraphics[width=0.3\textwidth]{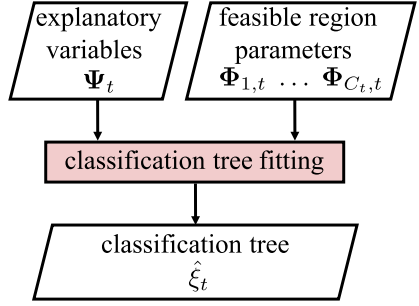} 
    \begin{correction_E}
    \caption{Illustration of the data and algorithm used to train $\hat \xi_t$, an approximation of $\xi_t$. The function $\hat \xi_t$ selects the best set of feasible region parameters given the expected values of a set of explanatory variables (\eg, outdoor temperature, day of the week). }  \label{fig:class_tree_alg}
        \end{correction_E}
\end{figure}

\correction{Recall that we use the load, indoor temperatures, and outdoor temperature relationships to group the elements of the training dataset. Naturally, we have no load nor indoor temperature data \emph{before} the new day.  However, each data point in the training set is associated with the explanatory variables denoted by $ \bd \Psi_t = \{ \bd \psi_{1,t}, \dots, \bd \psi_{K,t} \}$. The information encoded in $\bd \psi_{k,t}$ is anything that might influence the feasible region of the load during day $k$. For instance, $\bd \psi_{k,t}$ might include information on whether $k$ is a weekday, weekend, or a holiday, outdoor temperature during time $t$, solar irradiation levels, or building occupancy, among others. In this work, we use hourly outdoor temperatures, solar radiation, and day of the week as explanatory variables. However, considering a different set of explanatory variables might be appropriate in some cases and improve the effectiveness of the algorithm. }

We use expected values of the explanatory variables of day $K+1$ to select which $ \Phi_{c,t}$'s to use.  We assume there exists a function
\begin{equation*}
\xi_t: \bd \psi_{k,t} \rightarrow \left\{  \bd \Phi_{1,t}, \bd \Phi_{2,t}, \dots, \bd \Phi_{C_t,t}\right \}
\end{equation*} 
that maps the set of external data $\bd \psi_{k,t}$ to its associated set of parameters. Here, $\xi_t(\bd \psi_{k,t}) = \Phi_{c,t}$ when $k$ belongs in the set  $\mathcal{D}_{c,t}$ (recall that $\mathcal{D}_{c,t}$ is used to train $\Phi_{c,t}$). We estimate $\xi_t$ using a \emph{classification tree}~\cite{Loh}. Fig.~\ref{fig:class_tree_alg} 
 illustrates the training algorithm and inputs needed to train an approximation of $\xi_t$, which is denoted by $\hat \xi_t$.

Let $\left\{ \psi_{K+1,1}, \psi_{K+1,2},\hdots, \psi_{K+1,T}\right\}$ denote the set of explanatory variables for each time period of day $K+1$. Then, the predicted model to use is 
\begin{equation*}
\bd \Phi_{\hat c_t,t}= \hat \xi_t(\psi_{K+1,t}) 
\end{equation*} 
where $\hat \xi_t$ is a trained classification tree. The set of predicted feasible region parameters
\begin{equation*}
\left\{\bd \Phi_{\hat c_1,1}, \bd \Phi_{\hat c_2,2}, \dots, \bd \Phi_{\hat c_T,T}   \right\}
\end{equation*} 
describe the feasible region that models the building's flexibility during day $K+1$.

\begin{figure}
    \centering
    \includegraphics[width=0.485\textwidth]{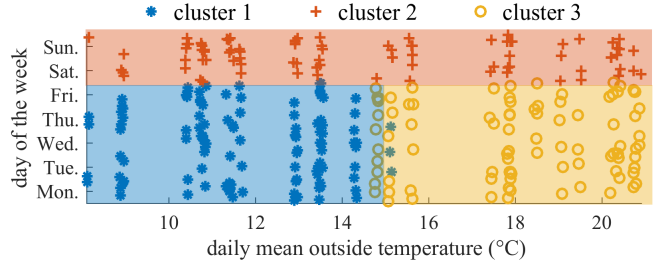} 
    \begin{correction_E}
    \caption{Training data grouped into three clusters. Each point in the plot represents a training data point. }\label{fig:choose_cluster}
    \end{correction_E}
\end{figure}
Take the data in Fig.~\ref{fig:choose_cluster} 
as an example.  In this case, the explanatory variables are the day of the week and mean outside temperature associated with each training data point. Notice that training data points in cluster $2$ come exclusively from weekends. On the other hand, training data points in clusters $1$ and $3$ come from colder and warmer weekdays, respectively. Then, a reasonable decision rule would be: use the parameters $\Phi_{ 2,t}$  day $K+1$ is a weekend; use $\Phi_{1,t}$ if $K+1$ is a weekday and the daily mean temperature is expected to be under $15^\circ$C; and use $\Phi_{ 3,t}$ otherwise. Fig.~\ref{fig:trained_tree} 
illustrated the decision tree structure that results from data in Fig.~\ref{fig:choose_cluster}
. We follow the same procedure for all $t$ from $1$ to $T$ and build the feasible region $\widehat{\bd{\mathcal{P}}}_{K+1}$ with the obtained sets of parameters. 

\begin{figure}
    \centering
    \includegraphics[width=0.45\textwidth]{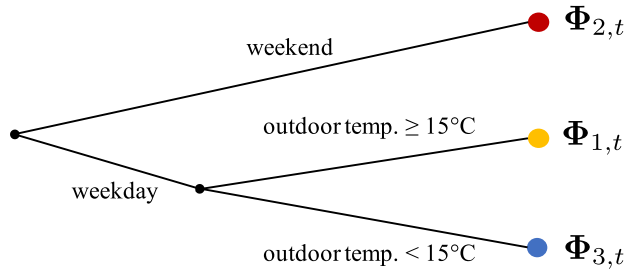} 
    \begin{correction_E}
    \caption{Decision tree used trained using points in Fig.~\ref{fig:choose_cluster} 
        .}  \label{fig:trained_tree}
    \end{correction_E}
\end{figure}

\end{correction_E}
	
	\section{Model Validation}
	\label{sec:case}

We test the proposed models using data obtained via EnergyPlus simulations of three different buildings from reference~\cite{CommercialRef}. Table~\ref{tab:title} provides a brief summary of important characteristics of each building type: peak load, average load, and thermal mass~(the amount of electric energy needed to cool the building by $1^\circ \mathrm{C}$).
\begin{table}[ht]
	\caption {Summary of building characteristics} \label{tab:title}
	\begin{center}
		\begin{tabular}{ |c|c|c|c| }
			\hline
			\textbf{Type} &   \textbf{Peak / avg. load} & \textbf{Thermal mass}    \\
			\hline
			Office 1 &  $27 /\ 10\; \mathrm{kW}$ &  $2.7\; \mathrm{ kWh /\ ^\circ C}$ \\
			\hline
			Office 2 & $15 /\ 6.7 \;\mathrm{kW}$ & $7.7\; \mathrm{ kWh /\ ^\circ C}$ \\
			\hline
			Supermarket & $140 /\ 86\;\mathrm{kW}$ & $77\; \mathrm{ kWh /\ ^\circ C}$ \\
			\hline
		\end{tabular}
	\end{center}
\end{table}

The sizes of the training, cross-validation, and test datasets for each building are $300$, $100$, and $100$, respectively. For the model selection stage, we use the day of the week (\eg, Monday), outdoor temperature, and solar irradiation as explanatory variables.

\subsection{Test error and optimal number of clusters}

\begin{figure}
		\centering
		\includegraphics[width=0.485\textwidth]{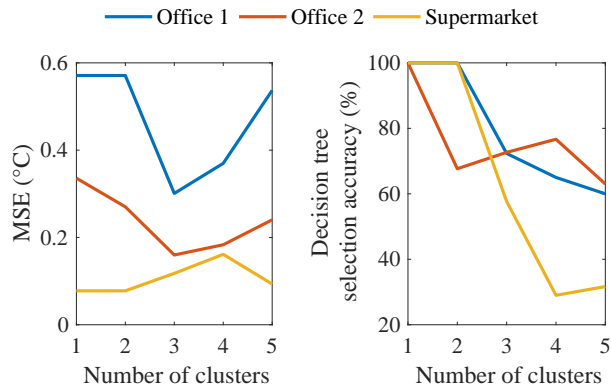} 
			\begin{correction_E}
		\caption{The left-hand plot shows root mean square error (for the cross-validation dataset) of the indoor temperature prediction for the three buildings at period $t=19$ as a function of the number of training data clusters. In this case, the optimal number of clusters is for office 1, office 2 and the supermarket are $3$, $3$, and $1$, respectively. The right-hand plot shows the accuracy of the decision tree as a function of the number of clusters. The accuracy tends to decline with the number of clusters. }  \label{fig:optimal_clusters}
	\end{correction_E}
\end{figure}

There is a trade-off when deciding the number of training data subsets $C_t$ to use. On the one hand, a small $C_t$ implies that more training data is available for each approximate feasible region.  On the other hand, a large $C_t$ implies that each approximate feasible region of the load models days that are more like each other.     \correction{Similar to references~\cite{Tibshirani_2005_Cluster, fu2017estimating} and as a special case of the hyperparameter tuning problem in machine learning, we define the optimal number of training data clusters $C^*_t$ as the number of clusters that minimizes the MSE of the temperature prediction functions over the cross-validation dataset. That is, the optimal number of clusters is the one that provides the best prediction over the cross-validation set. Then we use the test dataset (which is not used to tune the number of clusters) to measure the actual prediction performance.} For instance, the cross-validation error of temperature prediction for Office 1 is minimized when the number of training data clusters is $C^*_t=2$, as shown in the left-hand plot of Fig.~\ref{fig:optimal_clusters}
.  

The reason that the cross-validation initially decreases error with the number of clusters is that, as we divide the training data into more groups, each data cluster is used to approximate functions that are more like each other. However, at some point increasing the number of data clusters actually increases the cross-validation error (see the left-hand plot in Fig.~\ref{fig:optimal_clusters}). There are two main reasons for this phenomenon.  The first one is that a higher number of clusters means that each cluster contains fewer data points and thus the resulting estimator might be over-fitted. The second reason is that as the number of clusters increases, the accuracy of the decision tree $\hat \xi_t$ decreases (see the right-hand plot in Fig.~\ref{fig:optimal_clusters}) and the number of cross-validation data points predicted with the ``wrong" estimator increases. 

\begin{figure}
	\centering
	\includegraphics[width=0.485\textwidth]{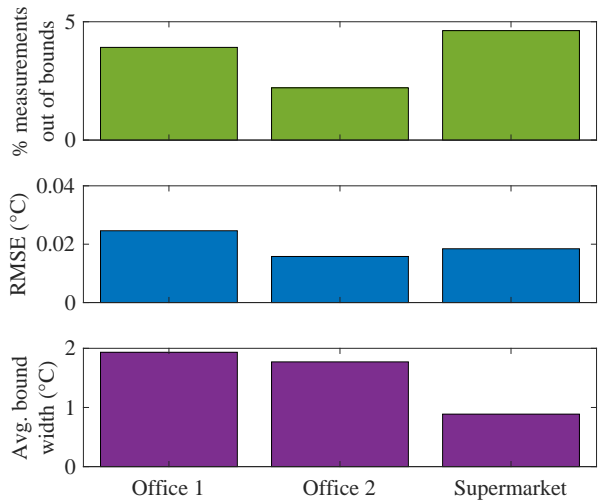} 
	\caption{Percentage of measurements out of bounds, test error of the temperature prediction, and average bound width for each building type. These statistics are computed using test set data.}  \label{fig:test_error}
\end{figure}

\subsection{Error analysis on the test set}

Fig.~\ref{fig:test_error} shows the test set percentage of indoor temperature measurements out of bounds, temperature prediction error, and tightness of the upper and lower predictions bounds. The percentage of indoor temperature measurements out of bounds is closely linked to the robustness parameter $\alpha$. Recall that an $\alpha$ level of robustness restricts the percentage of \emph{training} data points outside the prediction bands to be less than $100 \cdot \alpha \; \%$. As shown in Fig.~\ref{fig:test_error}, under the percentage of indoor temperature measurements out of bounds metric, the trained models perform well for the test set. 

The for the BLSE case, the RMSE metric is defined  by Eq.~\eqref{eq:RMSE_BLSE} in the Appendix. The errors are defined as zero if the temperature measurement is inside the prediction bands and as the distance to the nearest band if the measurement is not within the bounds (see Fig.~\ref{fig:BLSE_est}) for an illustration. The RMSE error for the proposed model is lower than the error given by the conventional RC circuit model (we compare our approach against the RC-circuit model in greater detail in Sec.~\ref{subsec:RC}).  However, the lower error achieved by our model is not entirely free. This lower error comes at the cost of more conservative modeling of the building's flexibility, as illustrated in Sec.~\ref{sec:wind_power_bal}. 

\begin{correction_E}
	\subsection{Comparison against the RC circuit model}    \label{subsec:RC}
	The most widely adopted alternative to the flexibility model offered in our work is the RC circuit model~\cite{Ma2012,Radecki,Gouda2002,Hao2015, Contreras-Ocana2016}. It expresses the indoor temperature change from time $t$ to $t+1$ a linear function of the indoor-outdoor temperature difference, HVAC power\footnote{In some cases, HVAC cooling/heating load is used in lieu of HVAC power, \eg~\cite{Ma2012}.}, and an independent thermal disturbance. Similar to~\cite{Ma2012,Contreras-Ocana2016},  the RC circuit model can be expressed as
	\begin{equation*}
	\phi_{t+1}^\mathrm{in}  - \phi_{t}^\mathrm{in} = A_t \cdot  (\phi_t^\mathrm{in} -\phi^\mathrm{out}_t) + B_t \cdot p_t^\mathrm{hvac} + D_t
	\end{equation*}
	where $A_t$ and $B_t$ relate indoor-outdoor temperature difference and HVAC power, respectively, to change in indoor temperature from $t$ to $t+1$. The independent thermal disturbance is denoted by $D_t$. For this case study, we estimate the $A_t$, $B_t$, and $D_t$ via a linear regression where the dependent variable is the temperature change $    \phi_{t+1}^\mathrm{in}  - \phi_{t}^\mathrm{in} $ and the regressors are $\phi_t^\mathrm{in} -\phi^\mathrm{out}_t$ and $p_t^\mathrm{hvac}$. 
	
	\begin{figure}
		\begin{correction_E}
			\centering
			\includegraphics[width=0.485\textwidth]{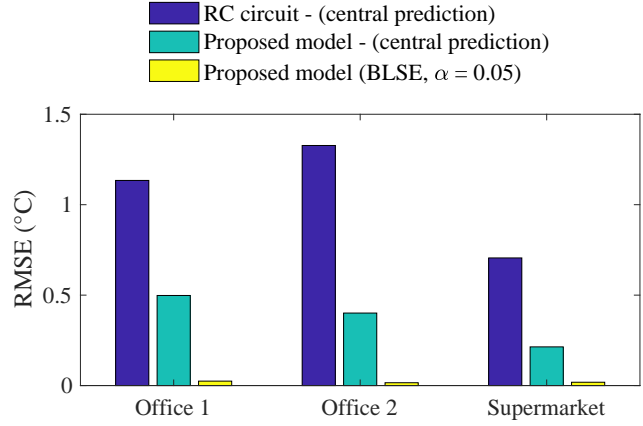} 
			\caption{Test set RMSE for three different models: the RC circuit model, the proposed model central prediction ($\alpha=1$) and the proposed model with $\alpha=0.05$. }  \label{fig:rc_comparison}
		\end{correction_E}
	\end{figure}
	
	Fig.~\ref{fig:rc_comparison} shows the test set RMSE for three different models: the RC circuit model, and our approach with $\alpha=1$  and $\alpha=0.05$. It is natural that the RMSE is orders of magnitude smaller with the $\alpha=0.05$ model: by definition, close to $95$\% of the predictions fall in the prediction band.  Less intuitive is the fact that our method with $\alpha=1$, equivalent to a central estimate, also outperforms the RC circuit model. There are two reasons for this. The first is that our use of clusters to model thermal dynamics with several linear functions instead of a single one. The second is that the traditional RC circuit model fails to use $p_{t}^\mathrm{hvac}$ to predict indoor temperature during time $t$~\cite{Ma2012,Radecki,Contreras-Ocana2016}.

\end{correction_E}

\section{Case study: building flexibility for wind power balancing}
\label{sec:wind_power_bal}
\begin{figure}
	\centering
	\includegraphics[width=0.485\textwidth]{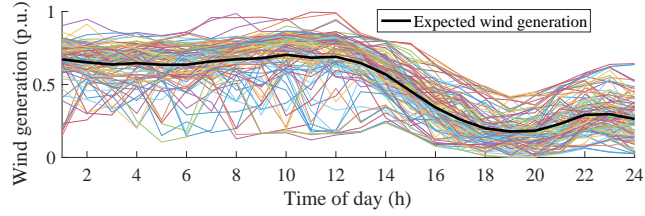} 
	\caption{Wind generation scenarios used in the case study. }  \label{fig:wind_scen}
\end{figure}
We consider a setting where an aggregator of buildings is contracted by a wind power producer to compensate deviations from the expected production. Let $\bd \Omega^\mathrm{w}$ represent the set of scenarios and each generation scenario be denoted by the vector $\bd \nu_{\omega^\mathrm{w}} \in \mathbb{R}_+^T$. The $t^\mathrm{th}$ entry of $\bd \nu_{\omega^\mathrm{w}}$ represents wind power at time $t$ in scenario $\omega^\mathrm{w}$. Then, the expected wind production is given by $\mathbb{E} [\bd \nu_{\omega^\mathrm{w}}]$ and the wind production deviation of scenario $\omega$ by $\bd \Delta_{\omega^\mathrm{w}} = \bd \nu_{\omega^\mathrm{w}} - \mathbb{E} [\bd \nu_{\omega^\mathrm{w}}]$. Fig.~\ref{fig:wind_scen} shows 100 wind scenarios obtained from references~\cite{Bukhsh2016, Pinson2013}.

In addition to wind uncertainty, we consider uncertainty in the building load. We model the stochastic component of building load $i$ using a $T$-dimensional normally distributed parameter $\bd \epsilon_i \sim \bd{\mathcal{N}}(\bd 0,\bd \Sigma_i)$.
Assuming that the stochastic components of the $N$ buildings are independent, the aggregate stochastic component of the load is $\sum_i^N \bd \epsilon_i = \bd \epsilon \sim \bd{\mathcal{N}}(\bd 0, \sum_i^N \bd \Sigma_i)$.
We represent $\bd \epsilon$ via scenarios $\{\bd \epsilon_{\omega^\mathrm{b}} \}_{\omega^\mathrm{b} \in \bd \Omega^\mathrm{b}}$.

	The aggregator's problem is as follows. In the first stage, \eg, in the day-ahead market, the aggregator schedules aggregate base load $\bd p^\mathrm{b} \in \mathbb{R}^T_+$ of the $N$ buildings at an energy price $\bd \tau \in \mathbb{R}^T$.  When the uncertainty in wind production materializes in the second stage, \eg, in the real-time, the aggregator can deviate from the base load to accommodate deviations and be remunerated by $v$ per unit energy. For instance, suppose that the wind production a particular hour is expected to be $10$ kWh but the actual production is $12$ kWh. To partially accommodate the $2$ kWh surplus, the building loads deviate from their base load of $50$ kWh to $51$ kWh. Then, the building pays $50\cdot \tau_t$ for day-ahead energy and receives $1 \cdot v$ for balancing services. The aggregator's problem can be written as follows: 
	\begin{subequations}\label{prob:case_study}
		\begin{align}
		&\min_{\substack{\bd p_i^\mathrm{b}, \bd p_{i,\omega^\mathrm{w}} \\ \bd p^\mathrm{b}, \bd p_{\omega^\mathrm{w}, \omega^\mathrm{b}} }}  \!\!\!\! \bd \tau^\top \bd p^\mathrm{b} + \mathbb{E}  [v\cdot| \bd p^\mathrm{b} - \bd p_{\omega^\mathrm{w}, \omega^\mathrm{b}} + \bd \Delta_{\omega^\mathrm{w}}|]  \label{eq:CSobj}
		\\
		&\mbox{s.t. }\nonumber \\
		& \bd p^\mathrm{b} = \sum_{i=1}^N \bd p^\mathrm{b}_i  \label{eq:Aggbase}
		\\
		& \bd p_{\omega^\mathrm{w}, \omega^\mathrm{b}}= \sum_{i=1}^N \bd p_{i,\omega^\mathrm{w}}  + \bd \epsilon_{\omega^\mathrm{b}} \;\forall\; \omega^\mathrm{w} \in \bd \Omega^\mathrm{w},\; \omega^\mathrm{b} \in \bd \Omega^\mathrm{b} \label{eq:Agg_w} 
		\\
		& \bd p^\mathrm{b}_i \in \widehat{\bd{\mathcal{P}}}_i \;\forall\; i=1,\hdots,N \label{eq:feas_region_1}
		\\
		&\bd p_{i,\omega^\mathrm{w}} \in \widehat{\bd{\mathcal{P}}}_i\;\forall\; i=1,\hdots,N, \; \omega^\mathrm{w} \in \bd \Omega^\mathrm{w}. \label{eq:feas_region_2}
		\end{align}
	\end{subequations}
	The objective function~\eqref{eq:CSobj}
	has two components: the cost of energy, $\bd \tau^\top \bd p^\mathrm{b}$, and the expected foregone revenue from balancing wind power deviations $\mathbb{E}  [v\cdot| \bd p^\mathrm{b} - \bd p_{\omega^\mathrm{w}, \omega^\mathrm{b}} + \bd \Delta_\omega|] $. The second stage variable $\bd p_{\omega^\mathrm{w}, \omega^\mathrm{b}}$ is the aggregate building load for wind scenario $\omega^\mathrm{w}$ and load uncertainty scenario $\omega^\mathrm{b}$. Eq.~\eqref{eq:Aggbase} defines the aggregate base load (first stage) as the sum of the base loads of each building. Similarly, Eq.~\eqref{eq:Agg_w}
	defines the aggregate load when scenarios $\omega^\mathrm{w}$ and $\omega^\mathrm{b}$ materialize (second stage) as the sum of individual loads. Finally, Eqs.~\eqref{eq:feas_region_1} and~\eqref{eq:feas_region_2}
	restrict the first and second stage load of each building, respectively, to be within their respective approximate feasible region. Recall that the feasible regions $\widehat{\bd{\mathcal{P}}}_i$ are defined by Eqs.~\eqref{equation:aprox_feasible_reg}
	and~\eqref{eq:u_l_estimate_functions} 
	as
		\begin{align*}
		\widehat{\bd{\mathcal{P}}}=\{  \bd{p} \; | \; & \hat{p}^\mathrm{min}_t \le   p_t \le \hat{ p}_t^\mathrm{max}  \; \forall\; t=1, \dots, T \\
		& \hat{ \theta}^\mathrm{U}_{t}(\phi_0^\mathrm{in}, \phi_t^\mathrm{out}, \bd p_{1:t})  \le \hat{\theta}_{t}^\mathrm{max} \; \forall \; t=1,\dots,T \\
		&\hat{\theta}^\mathrm{L}_{t}(\phi_0^\mathrm{in}, \phi_t^\mathrm{out}, \bd p_{1:t})  \ge \hat{\theta}_{t}^\mathrm{min}\; \forall \; t=1,\dots,T \\
		& \bd{p} = \begin{bmatrix} p_1&\hdots&p_T\end{bmatrix}^\top,\; \bd p = \bd p_{1:T} \}
		\end{align*} 
	and its parameters  are determined using the  estimation procedures outlined in Section~\ref{sec:data} and detailed in the Appendix.

Problem~\eqref{prob:case_study} is formulated as a stochastic linear program and modeled using Julia's JuMP environment~\cite{DunningHuchetteLubin2017}. The problem is solved using Gurobi Optimizer~\cite{gurobi} on a desktop computer running on a Intel(R) Xenon(R) CPU E3-1220 v3 @ 3.10 GHz with 16 GB of RAM.

\begin{correction_E}
	\subsection{Wind forecast error mitigation and balancing compensation}
	
	\begin{figure}
		\centering
		\includegraphics[width=0.485\textwidth]{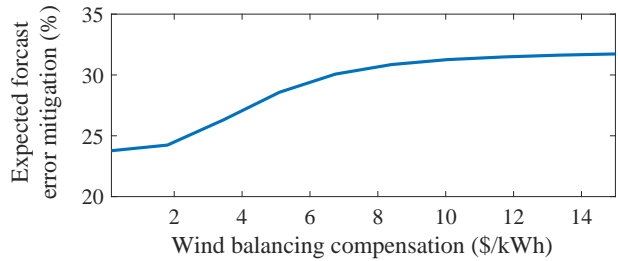}
		\begin{correction_E}
		\caption{Expected forecast error mitigation by all 3 buildings as a function of compensation, $v$. }  \label{fig:case_study2}
	\end{correction_E} 
	\end{figure}
	\begin{figure}
		\centering
		\includegraphics[width=0.485\textwidth]{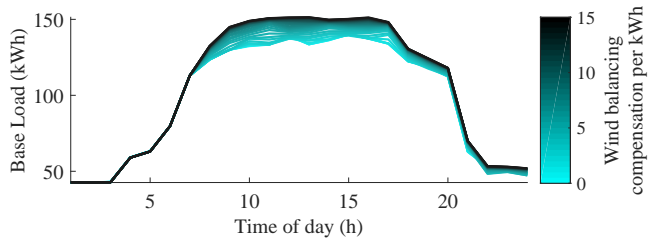} 
		\begin{correction_E}
		\caption{Base load at different wind balancing compensation levels. } \label{fig:case_study3}
	\end{correction_E} 
	\end{figure}
	
	Let the cost of energy be $1$ throughout the day and the installed wind capacity be one-third of the peak load. Depending on the compensation for wind balancing, the three buildings can mitigate around $25$--$30$ \% of the wind forecast errors.  As expected, and as shown in Fig.~\ref{fig:case_study2}, the amount of mitigated forecast error increases as the compensation for wind balancing increases. This result can be explained as follows. When the compensation is low, the base load tends to be low in order to minimize energy costs (see the lighter shades in Fig.~\ref{fig:case_study3}.  In this case, the low base load is poorly positioned to be further decreased in real-time to compensate wind shortages. As the balancing compensation increases, however, it becomes economically attractive for the building to position its base load at higher levels and increase its ability to accommodate wind shortages. The effect of the balancing compensation on the base load is shown in Fig.~\ref{fig:case_study3}.
	
\end{correction_E}

\begin{correction_E}
	\subsection{Demonstration of robustness and tractability}
	
	Robustness and tractability of are the two central characteristics of our model. The former claims that a building load profile deemed feasible by our model will not violate temperature limits during the actual building operation (to a degree of confidence determined by $\alpha$).   The latter claims that our model can be easily, and without significant computational burden, be incorporated into typical power system analysis frameworks (such as the one presented in this case study). 
	
	\begin{figure}
		\begin{correction_E}
			\centering
			\includegraphics[width=0.485\textwidth]{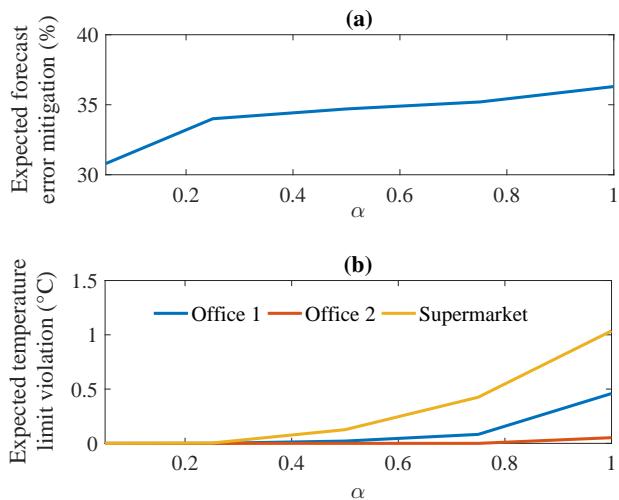}
			\caption{Plot a shows the expected forecast error mitigation by all 3 buildings and expected temperature violation by each building as a function of the robustness parameter $\alpha$. Plot b shows the expected temperature limit violation of each building as a function of $\alpha$. Notice that there is a trade-off between error mitigation and robustness to indoor temperature prediction errors. }  \label{fig:case_study_alpha}
		\end{correction_E}
	\end{figure}

	First, we analyze the effect of the robustness parameter $\alpha$ on the operation of the building load. Recall that a small $\alpha$ produces a more robust model (fewer measurements fall outside the prediction band) and a large $\alpha$ produces a less robust model. As shown in Fig.~\ref{fig:case_study_alpha}(a), the expected forecast error mitigation increases with $\alpha$. That is, as the robustness of the model decreases, it allows more aggressive operation of the building load to compensate forecast errors.  However, less robust models such as the RC circuit model, risk allowing load profiles that are not feasible during operation of the building (see an illustration of this phenomenon in Fig.~\ref{fig:model_vs_RC_illust}). As shown in Fig.~\ref{fig:case_study_alpha}(b), as the robustness parameters $\alpha$ increases so does the expected indoor temperature limits violations.  All in all, the user faces a trade-off when tuning the robustness parameter: a larger $\alpha$ allows for more aggressive operation of the HVAC load for forecast error mitigation but also poses a higher risk of causing indoor temperature limit violations.

	\begin{figure}
		\begin{correction_E}
			\centering
			\includegraphics[width=0.485\textwidth]{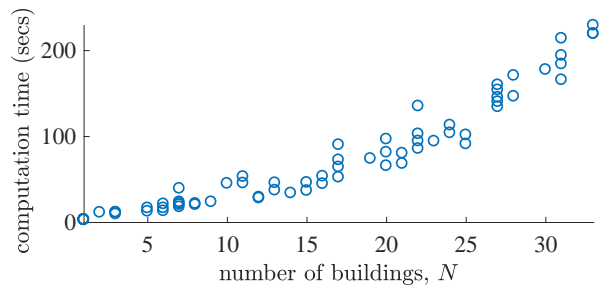}
			\caption{ Time required to solve Problem~\eqref{prob:case_study} 
				under different number of buildings. Note that each building is represented by 101 scenarios : one base load, $ \bd p^\mathrm{b}_i$, and one for each wind power scenario, $\bd p_{i,\omega^\mathrm{w}} $.}  \label{fig:comp_time}
		\end{correction_E}
	\end{figure}

	We demonstrate that the proposed model is tractable by increasing the number of buildings in Problem~\eqref{prob:case_study} 
	and showing that the computational burden to solve it remains manageable (see Fig.~\ref{fig:comp_time}). It is worth noting that in this case study, each building is represented by 101 scenarios: one base load, $ \bd p^\mathrm{b}_i$, and one for each wind power scenario, $\bd p_{i,\omega^\mathrm{w}}$. Thus, Problem~\eqref{prob:case_study} 
	with $N=10$ building case, for instance, is equivalent to solving a deterministic problem that involves $1010$ buildings.
	
\end{correction_E}

	\section{Conclusion}
	\label{sec:conclsion}
	
	In this work, we propose a method to estimate the feasible region of the building load using simple linear relations that is robust to temperature prediction errors. Our method ensures that a building's HVAC system is able to maintain acceptable occupant comfort while providing flexibility to the power system. Its mathematical simplicity makes our model tractable in the sense that it can be easily incorporated into optimization and control environments that are common in power systems applications. For instance, the proposed model can be seamlessly incorporated into problems such as the economic dispatch, demand response scheduling and control, optimal power flow,  and unit commitment, among others. The proposed model requires relatively little data to be trained and can be applied to all types of buildings, without requiring detailed sensor data. We compare our model to the RC circuit model and demonstrate a practical application in which an aggregator uses three different buildings (two offices and a supermarket) to balance wind generation forecast errors.

	\appendices
\section{Clustering the training data} \label{app:data}
Let the the vector $\bd w_{k,t}$ be defined $\forall\; k\in \bd{\mathcal K}$ and $t=1,\hdots,T$ as $\bd w_{k,t} = \begin{bmatrix}  \bd p_{k,1:t}^\top & \phi_{0,t}^\mathrm{in}   &\phi_{k,t}^\mathrm{in}&\phi_{k,t}^\mathrm{out} \end{bmatrix}^\top$.  Since the vectors $ \bd w_{k,t}$ contain data on different units and potentially different magnitudes, normalizing the data prevents the clustering algorithm from unfairly assigning more importance to some of the elements of $\bd w_{k,t}$.  Denote a normalized matrix of horizontal concatenation of all $\bd w_{k,t}$'s as $\bd W_t=\mathrm{norm}( \begin{bmatrix} \bd w_{1,t}& \bd w_{2,t}& \hdots & \bd w_{|\mathcal K|,t} \end{bmatrix} )$. The matrix $\bd W_t$ is normalized such that the mean of each row is zero and the $\ell_2$ norm of each row is $1$. We use the K-means algorithm~\cite{macqueen1967} to group the columns of $\bd W_t$ matrix into $C_t$ separate clusters.  The indices of $\bd w_{k,t}$'s assigned to cluster $c$ are denoted by  $\bd{\mathcal{K}}_{c,t}$.

\section{Estimate of the temperature and load limits} \label{app:est}
We estimate  $ \theta_{c,t}^\mathrm{max}$ as the maximum indoor temperature during time period $t$ during days in the set $\bd{\mathcal{K}}_{c,t}$, \ie, $\hat{\theta}_{c,t}^\mathrm{max} = \max(\{ \phi^\mathrm{in}_{k,t}\}_{k \in \bd{\mathcal{K}}_{c,t}})$.
The minimum temperature limit, the upper and lower load bounds for each cluster are estimated using an analogous procedure, \ie, $\hat{\theta}_{c,t}^\mathrm{min} = \min(\{ \phi^\mathrm{in}_{k,t}\}_{k \in \bd{\mathcal{K}}_{c,t}})$, $\hat{p}_{c,t}^\mathrm{min} = \min(\{ p_{k,t}\}_{k \in \bd{\mathcal{K}}_{c,t}})$, $\hat{ p}_{c,t}^\mathrm{max} = \max(\{ p_{k,t}\}_{k \in \bd{\mathcal{K}}_{c,t}})$.

		\section{Bounded least squares estimation}\label{app:blse}
		Let an estimate of upper bound of the indoor temperature at time $t$ and day $k$ be an affine function of the initial temperature $\phi_{k,0}$, the outdoor temperature at time $t$, $\phi_{k,t}^\mathrm{out}$, and load form the first to the $t^\mathrm{th}$ time period,
		\begin{align*}
		\hat{\theta}_{k,t}^\mathrm{U}(\phi_{k,0}^\mathrm{in}, \phi^\mathrm{out}_{k,t}, \bd{p}_{k, 1:t})  = \overline{\bd{a}}_{c,t}^{ \top}\bd{p}_{k, 1:t}+\overline{\bd{b}}_{c,t}^{\top}\begin{bmatrix}\phi_{k,0}^\mathrm{in} & \phi^\mathrm{out}_{k,t}  & 1\end{bmatrix}^\top.
		\end{align*}
		Similarly, the lower bound estimate of the indoor temperature at time $t$ is
		\begin{align*}
		\hat{\theta}_{k,t}^\mathrm{L} (\phi_{k,0}^\mathrm{in}, \phi^\mathrm{out}_{k,t}, \bd{p}_{k, 1:t}) = \underline{\bd{a}}_{c,t}^{\top}\bd{p}_{k, 1:t}+\underline{\bd{b}}_{c,t}^{\top}\begin{bmatrix}\phi_{k,0}^\mathrm{in} & \phi^\mathrm{out}_{k,t}  & 1\end{bmatrix}^\top.
		\end{align*}

		We cast the problem of finding values of $\overline{\bd a}_{c,t}$, $\underline{\bd a}_{c,t}$, $\overline{\bd b}_{c,t}$, and $\underline{\bd b}_{c,t}$ such that the square error and a measure of the tightness of the bounds are minimized as the following convex quadratic program:
		\begin{subequations} \label{eq:BLSE_prob}
			\begin{align}
			&\argmin_{ \substack{\overline{\bd a}_{c,t}, \underline{\bd a}_{c,t}, \overline{\bd b}_{c,t}, \underline{\bd b}_{c,t} \\ \hat{\theta}_{k,t}^\mathrm{U}, \hat{\theta}_{k,t}^\mathrm{L}  \\ J_{k,t}^\mathrm{U} ,J_{k,t}^\mathrm{L}, J_t^\mathrm{A} } } \beta_i \cdot\sum_{k\in \bd{\mathcal{K}}_c} \left(J_{k,t}^{\mathrm{U}} +J_{k,t}^{\mathrm{L}} \right)^2 + (1 - \beta_i)\cdot J^\mathrm{A}_t \\
			& \mbox{s.t.} \nonumber \\
			& \hat{\theta}_{k,t}^\mathrm{U} = \overline{\bd{a}}_{c,t}^{\top}\bd{p}_{k, 1:t}+\overline{\bd{b}}_{c,t}^{\top}\begin{bmatrix}\phi_{k,0}^\mathrm{in} & \phi^\mathrm{out}_{k,t}  & 1\end{bmatrix}^\top  \forall\; k \in \bd{\mathcal{K}}_{c,t} 
			\label{eq:BLSE_ub}
			\\
			& \hat{\theta}_{k,t}^\mathrm{L} = \underline{\bd{a}}_{c,t}^{\top}\bd{p}_{k, 1:t}+\underline{\bd{b}}_{c,t}^{\top}\begin{bmatrix}\phi_{k,0}^\mathrm{in} & \phi^\mathrm{out}_{k,t}  & 1\end{bmatrix}^\top \; \forall\; k \in \bd{\mathcal{K}}_{c,t}
			\label{eq:BLSE_lb} 
			\\
			& J_{k,t}^\mathrm{U}  \ge \phi^\mathrm{in}_{k,t} - \hat{\theta}_{k,t}^\mathrm{U}  \; \forall \; k \in\bd{\mathcal{K}}_{c,t} \label{eq:BLSE_U_error}
			\\
			& J_{k,t}^\mathrm{L}  \ge   \hat{\theta}_{k,t}^\mathrm{L} - \phi^\mathrm{in}_{k,t}\; \forall \; k \in \bd{\mathcal{K}}_{c,t}
			\label{eq:BLSE_L_error}
			\\
			&J^\mathrm{A}_t =  \sum_{k \in \bd{\mathcal{K}}_{c,t}}\hat{\theta}_{k,t}^\mathrm{U} - \hat{\theta}_{k,t}^\mathrm{L} 
			\label{eq:BLSE_tb} 
			\\
			&\hat{\theta}_{k,t}^\mathrm{U} \ge \hat{\theta}_{k,t}^\mathrm{L} \; \forall\; k \in \bd{\mathcal{K}}_{c,t} \label{eq:BLSE_u>l}
			\\
			& \overline{\bd a}_{c,t} \le 0,\; \underline{\bd a}_{c,t} \le 0. \label{eq:BLSE_neg_param}
			\\ 
			& J_{k,t}^\mathrm{U} \ge 0  \; \forall\; k \in\bd{\mathcal{K}}_{c,t} \\
			&J_{k,t}^\mathrm{L} \ge 0  \; \forall\; k \in \bd{\mathcal{K}}_{c,t} \label{eq:BLSE_error_pos} 
			\end{align}
		\end{subequations}
		The objective function of problem~\eqref{eq:BLSE_prob} is composed of two weighted components: \emph{1)} the sum of squared errors and \emph{2)} a measure of the tightness of the upper and lower estimates. The first component is weighted by $\beta_i$ while the second one is weighted by $1-\beta_i$ where $\beta_i \in (0,1)$.  When as $\beta_i \rightarrow 1$, the bounds become wider and more points fall within them. Conversely, when $\beta_i$ is small, the bounds are tighter.
		
		Equations~\eqref{eq:BLSE_ub} and~\eqref{eq:BLSE_lb}, define the estimates of the upper and lower estimates, respectively. Eq.~\eqref{eq:BLSE_U_error} 
		defines the upper estimate error $J_{k,t}^\mathrm{U}$ to be the distance between the upper estimation and the measurement if this quantity is positive and zero otherwise. Similarly, Eq.~\eqref{eq:BLSE_L_error} 
		defines the lower bound error $J_{k,t}^\mathrm{L}$ to be the distance between the measurement and the lower bound estimation if this quantity is positive and zero otherwise.
		
		The measure of the tightness of the prediction band $J_t^\mathrm{A}$ is defined in Eq.~\eqref{eq:BLSE_tb} 
		as the sum of the distance between the lower and upper estimates over all samples $\bd{\mathcal{K}}_{c,t}$. We restrict the upper estimate to be higher than the lower bound in Eq.~\eqref{eq:BLSE_u>l}
		. Without loss of generality, we assume that every training day is either cooling day. Then, everything else equal, higher load must translate into lower temperature. Therefore  $\overline{\bd a}_{c,t}$ and $\underline{\bd a}_{c,t}$ are restricted to be negative as in Eq.~\eqref{eq:BLSE_neg_param}
		. If the training days are all heating days, the signs in Eq.~\eqref{eq:BLSE_neg_param}
		are reversed. Finally, we define the root mean square error of the BLSE as 
		\begin{equation} \label{eq:RMSE_BLSE}
		\mathrm{RMSE}_{c,t} = \sqrt{\frac{\sum_{k \in \bd{\mathcal{K}}_{c,t}} \left(J^\mathrm{U}_{k,t} +J^\mathrm{L}_{k,t}  \right)^2}{|\bd{\mathcal{K}}_{c,t}|}}.
		\end{equation}

\subsection*{The BLSE algorithm}
Let $\mathrm{BLSEF}(\beta_i)$ denote a function that takes the scalar $\beta_i \in (0,1)$ and solves problem~\eqref{eq:BLSE_prob} 
and outputs the optimal values of  $\overline{\bd{a}}_{c,t}$ $\underline{\bd a}_{c,t,i}$, $ \overline{\bd{b}}_{c,t,i}$, $\underline{\bd b}_{c,t,i}$, $J^\mathrm{A}_{t,i}$, and calculates the percentage of training measurements that are higher than the upper estimate or lower than the lower estimate  $\pi^\mathrm{out}_{t,i}$. The percentage of out of prediction band measurements is calculated as
\begin{equation*}\pi^\mathrm{out}_{t,i} = \frac{\sum_{k \in \bd{\mathcal{K}}_{c,t}}\mathbb{I}(\hat{\theta}^\mathrm{U}_{k,t} \le \theta_{k,t} \mathrm{\;or\;} \hat{\theta}^\mathrm{L}_{k,t} \ge \theta_{k,t})}{|\bd{\mathcal{K}}_{c,t}|}\end{equation*} 
where $\mathbb{I}(\cdot)$ is the indicator function.

Now we describe the BLSE algorithm (see Algorithm~\ref{alg:blse} below). Its inputs are: the training data set, a maximum out of band percentage $\alpha$ (\eg, $5 \%$), and a vector $\bd \beta = \{{0, \frac{1}{M-1},\frac{2}{M-1},\hdots, 1} \}$.  The parameter $M$ is an integer greater than $1$ to be selected by the modeler\footnote{A small $M$ reduces the computation time of Algorithm~\ref{alg:blse}
	but might yield less accurate solutions. A large $M$, on the other hand, increases the computation time but yields a more accurate solutions. In this work we use $M=100$}.  The outputs of the BLSE algorithm are the trained parameters of $\hat{\theta}_{k,t}^\mathrm{L}$ and $\hat{\theta}_{k,t}^\mathrm{U}$. For each time $t$ it does the following: it goes through each element of $\bd \beta$, $\beta_i$, it solves $\mathrm{BLSEF(\beta_i)}$. Then, among the solutions that yield an out-of-band percentage smaller than $\alpha$,  it selects the one that solution that yields tighter band, \ie, the smallest $J^\mathrm{A}_{t,i}$. 


\begin{algorithm}
	
	\KwIn{ $\{\bd{\phi}_k^\mathrm{in}, \bd{ p }_k , \bd{\phi}_k^\mathrm{out}, \phi_{0,k}^\mathrm{in} \}_{k \in \bd{\mathcal{K}}}$, $\alpha$ , $\bd \beta = \{{0, \frac{1}{M-1},\frac{2}{M-1},\dots, 1}  \},  \{\bd{\mathcal{K}}_{1,t}, \dots, \bd{\mathcal{K}}_{C_t,t} \}_{t=1,\dots, T} $ }
	\KwOut{  $\{\{\overline{\bd{a}}_{c,t}^*, \underline{\bd{a}}_{c,t}^*, \overline{\bd{b}}_{c,t}^*, \underline{\bd{b}}_{c,t}^*\}_{c=1\dots, C_t} \}_{t=1,\dots, T}$ }
	\For{$t=\{ 1,\dots, T\}$}{
		\For{$i=\{ 1,\dots, M\}$}{
			\For{$c=1\cdots, C_t$}{ 
				$(\overline{\bd{a}}_{c,t,i},\;\underline{\bd{a}}_{c,t,i}, \overline{\bd{b}}_{c,t,i},\;\underline{\bd{b}}_{c,t,i},J^\mathrm{A}_{t,i}, \pi^\mathrm{out}_{t,i} )= \mathrm{BLSEF}(\beta_i)$ \\
				
			}
			$i^*=\argmin_{\pi^\mathrm{out}_{i,t} \le \alpha}\{J^\mathrm{A}_{t,i}  \}_{i=1,\dots,M} $ \\
			
			$(\overline{\bd{a}}_{c,t}^*,\;\underline{\bd{a}}_{c,t}^*, \overline{\bd{b}}_{c,t}^*,\;\underline{\bd{b}}_{c,t}^*) =(\overline{\bd{a}}_{c,t,i^*},\;\underline{\bd{a}}_{c,t,i^*}, \overline{\bd{b}}_{c,t,i^*},\;\underline{\bd{b}}_{c,t,i^*})$
		}
		
	}
	\caption{The bounded least square estimation (BLSE) algorithm.}
	\label{alg:blse}
	
\end{algorithm}

Note that a larger $M$ will increase the computation time required to run Algorithm~\ref{alg:blse} 
but will produce a larger set $\{\pi^\mathrm{out}_{i,t}\}_{1=1,\dots,M}$. A larger set of  $\pi^\mathrm{out}_{i,t}$'s makes it likelier that the optimal $\pi^\mathrm{out}_{i^*,t}$ is closer to the desired robustness parameter $\alpha$.

	\bibliographystyle{IEEEtran}

	\bibliography{bibliography}

\end{document}